\newcommand{\beq}{\begin{quote}}
\newcommand{\enq}{\end{quote}}
\newcommand{\be}{\begin{equation}}
\newcommand{\en}{\end{equation}}
\newcommand{\del}{\delta}
\newcommand{\om}{\omega}
\newcommand{\Om}{\Omega}

\newcommand{\eps}{\epsilon}
\documentstyle[12pt,epsfig]{article}
\begin{document}
\title{ Continuity and Stability of families of figure eight orbits 
with finite angular momentum. 
} 
\date{}
\author{Michael Nauenberg\\
Department of Physics\\
University of California, Santa Cruz, CA 95064 
}
\maketitle
\begin{abstract}

Numerical solutions are presented for a family
of  three dimensional periodic orbits with three equal masses
which  connects the classical circular orbit of Lagrange with 
the recently discovered planar figure eight orbit 
with zero total angular momentum \cite{moore},\cite{mont1}. 
Each member of  this  family is an orbit with  finite angular
momentum that is  periodic  in a frame 
which rotates with frequency
$\Omega$ around the horizontal symmetry axis of the
figure eight orbit. Numerical solutions  for figure eight shaped
orbits with finite angular momentum orbits 
were first reported in \cite{michael1}, and
mathematical proofs for the existence of such orbits
were given in \cite{march2}, and more recently in \cite{mont2} 
where also some numerical solutions have been presented. 
Numerical evidence is given here that the family of such 
orbits is a continuous function of the rotation frequency $\Omega$ 
which varies between $\Omega=0$ ( for the planar figure eight 
orbit with intrinsic frequency $\omega$), 
and $\Omega =\omega$ (for a circular Lagrange orbit) 
Similar numerical solutions are 
also found for $n>3$ equal masses, where $n$ is an odd integer,
and an illustration is given for $n=21$.
Finite angular momentum  orbits
were also obtained numerically
for rotations along the two other symmetry axis of the
figure eight orbit \cite{michael1}, and some new results are
given here. Recently, existence proofs for families of 
such orbit and further 
numerical solutions have been given in \cite{mont2}.
The stability of these orbits is examined numerically 
without the restriction to a linear approximation \cite{mont2}, 
and some examples are given  of nearby stable orbits 
which bifurcate from these families.

\end{abstract}

\subsection*{Introduction}

In 1993 C. Moore \cite{moore} 
discovered a new periodic orbit  for three particles
of equal mass  moving synchronously
on a symmetric planar figure eight orbit under the action of their mutual
attractive gravitational forces.
Subsequently, a mathematical proof for the existence of this orbit and 
an accurate numerical solution were  given in \cite{mont1}.
This orbit  has total angular momentum equal to zero.  Similarly 
shaped figure eight orbits with finite angular momentum were 
obtained numerically, and  were first presented in  \cite{michael1}.
Before  Moore's discovery, C. Marchal\cite{march1} 
had  obtained a perturbative orbit in the shape
of a three dimensional figure eight which bifurcates 
from the classical Lagrange orbit for three equal masses 
located on the vertices of an equilateral triangle which rotates uniformly 
around its center. This  perturbed orbit is  periodic in a 
frame rotating  around an axis normal 
to the plane and through the center of the Lagrange orbit, 
with a small rotation angular frequency proportional 
to the square of the maximum
amplitude of the perturbation amplitude along this axis.
After  learning about the existence of the 
planar figure eight orbit, Marchal gave a  proof 
\cite{march2} for  the 
existence of a family of three dimensional orbits where
each member is periodic in a frame rotating 
around the horizontal symmetry axis of the planar figure
eight orbit. This family connects the figure eight orbit to a
Lagrange orbit, with its plane normal to this axis, by varying 
the angular rotation frequency $\Omega$ of this frame.
Both the perturbative orbit that bifurcates from the
Lagrange orbit \cite{march1},\cite{march2}, and the finite angular momentum
orbit that  bifurcates from the figure eight orbit
along the horizontal axis of symmetry \cite{michael1},  
are members of this family.
Numerical solutions of periodic orbits 
that bifurcate from the  figure eight
in a frame rotating around an axis  
normal to its horizontal symmetry axis, 
were also reported in \cite{michael1}. 
Recently, a proof for the existence of 
of such periodic orbits  
and some numerically  generated figures for these orbits
were given in reference \cite {mont2}.
But a proof has not been given for the {\it continuity} as a function of the
rotation frequency $\Omega$ of the  family of periodic  orbits  
between the Lagrange orbit and 
the planar figure eight orbit. In reference \cite{march2},  Marchal 
concluded  that '' the next step will be of course...
the numerical verification of the continuity up to the
eight shaped orbit... ''. Moreover,  the stability of these 
orbits have been examined only in a linear approximation \cite{mont2},
but this  turns out to be inadequate, as will be shown here. 

In this paper, we extend previous numerical results 
\cite {michael1}, \cite{mont2}, and  
show that the family of three dimensional 
orbits which  connects the classical circular orbit of Lagrange with 
the planar figure eight orbit
is a continuous function of the  angular rotation frequency 
$\Omega$ of a frame rotating around
the horizontal symmetry axis of the figure eight,
where the orbits are periodic.   
These orbits are obtained by expanding its coordinates 
in the  rotating frame in
a Fourier series, and calculating the Fourier  coefficients
by a  steepest descent gradient method \cite{michael1}
discussed in section I. 
The first nine  coefficients 
for this expansion are given in Figs. 4-6, 
which demonstrates
that as a function of the rotation frequency $\Omega$ of this frame (with
fixed  frequency $\omega=1$ of the periodic orbit), 
the coefficients vary continuously between those for the planar 
figure eight orbit for $\Omega=0$, and those  for the Lagrange orbit 
at  $\Omega=1 $. In particular,
we verify numerically that in these two limits, the
Fourier  coefficients approach the values obtained from 
perturbative solutions. Marchal evaluated the first
few orders  for the  perturbative solution  near the Lagrange
orbit, which can be expressed  as a series in
powers of the {\it square root} of
the rotation frequency $\Omega$  \cite {march1},\cite {march2}.  
The corresponding perturbation solution
near the zero angular momentum figure eight orbit
is in powers of $\Omega$, 
as was pointed out in \cite{michael1}.
In section II we illustrate the application of the
Fourier gradient method near the classical Lagrange solution 
by a first order perturbation calculation similar to Marchal's.
As is well known, this solution 
can be extended to any odd number $n$ of equal  masses 
on the vertices of a  regular $n${th} polygon,
and this is also the case for the figure eight 
shaped orbit, as was  shown numerically 
in \cite{simo1}, \cite{michael1}.
Likewise, we have found that the family of  
periodic orbits connecting these two orbits 
can also be extended to an odd number of masses.
As an illustration, we present in Fig. 2 numerical results
for the case that  $n=21$.
Some additional  numerical solutions for
finite angular momentum orbits
associated with rotations of  the planar figure eight orbit 
along the $y$ and $z$ axes of symmetry  \cite{michael1}
\cite{mont2} are shown in Figs. 7-9. 
These results are discussed in section III. Finally, in section IV 
we examine the stability of these families of orbits by integrating
the equations of motion numerically with a fourth order
Runge-Lenz algorithm, taking for initial conditions values 
of the coordinates and momenta nearby to
those obtained previously with our Fourier gradient method 
for periodic solutions.

\subsection* {I. Fourier gradient method}

The symmetry properties of the periodic planar figure eight shaped  orbit
for three equal masses \cite{moore},\cite{mont1},\cite{march2} 
dictate the form of the Fourier series expansion of the orbital coordinates.
For Newtonian gravitational interactions, the coordinate
and the time scale according to the Keplerian  transformation 
$x(t)= \tau ^{2/3}f(t/\tau)$, 
where $\tau$ is the period;  and similarly for the 
$y(t)$ coordinate. Choosing  for convenience  $\tau=2\pi$, 
we have   
\be
\label{xt}
x(t)=\sum_{k=1}^{\infty} a(k_o)sin(k_ot),
\en
\be
\label{yt}
y(t)=\sum_{k=1}^{\infty} b(k_e)sin(k_et)
\en
where $k_o=2k-1$  and $k_e=2k$.
This represention is also valid for any odd number
$n\geq 3$ of equal masses, where the position of the $j${\it th} mass at time
$t$  is given by $x_j(t)=x(t+j\tau/n)$
and $y_j(t)=y(t+j\tau/n)$ for $j=0,..(n-1)$.
The conservation of total linear momentum requires that 
the Fourier coefficients $a(k_o)$  and $b(k_e)$ vanish
for $k_o$ mod $n$ and $k_e$ mod $2n$ respectively.
For a figure eight shaped  orbit with finite
angular momentum about the $x$ or $y$ axis of symmetry, 
this expansion remains valid in a frame rotating around
one of these symmetry axes,
but then the orbit depends also  on the  coordinate $z$ 
normal to the $x-y$ plane in the rotating frame.  
For the case of rotations around the $x$ axis, 
it can be readily shown from the equations of motion
that the Fourier expansion of $z(t)$  has the form
\be 
\label{zte}
z(t)=\sum_k c(k_e)cos(k_et),
\en
while for rotations around the $y$ axis
\be
\label{zto}
z(t)=\sum_k c(k_o)cos(k_ot).
\en
In the fixed frame, the corresponding  coordinates
can  conveniently  be
obtained by introducing complex variables. For example,
for rotation around the $x$ axis,  we have
\be
x_j(t)=\sum_k sin(k_o(t+2\pi j/n)),
\en
and
\be
\label{complex1}
y_j(t)+iz_j(t)= e^{i\Omega t} \sum_k e^{ik_e (t+2\pi j/n)} (b(k_e)+ic(k_e)), 
\en
where $\Omega$ is the angular rotation frequency of the frame.
Setting $\Omega=1-\beta$ for $0 \leq \Omega$  and  $\Omega=-1+\beta$
for $\Omega \leq 0$, we see that an orbit which is periodic
in the coordinate rotating with angular frequency $\Omega$
is also periodic in another frame rotating with angular frequency $\beta$.
For $n=3$ and $\Omega=-1+\beta$, this representation
corresponds to that given by Marchal \cite{march1},\cite{march2}, 
where $\beta$ is the angular rotation frequency of a frame relative to
which the Lagrange orbit has period $2\pi$,  

We evaluate these Fourier coefficients by  
determining the extrema of the action integral $A$ 
for an orbit given by Eqs. \ref{xt}-\ref{zto} in a 
frame rotating with angular frequency $\Omega$. For example,
for  the case that the axis of rotation is taken along the $x$ axis
of symmetry of the planar figure eight orbit, we have 
\be
\label{action1}
A=\int_0^{2\pi} dt [K(t)-P(t)+ \frac{1}{2} \Om^2 I_x(t)+\Om L_x(t)],
\en
where  $K(t)$ is the kinetic energy, $P(t)$ the potential energy,
$I_x(t)$  the $x$-moment of inertia, and $L_x(\theta)$  the $x$ component 
of angular momentum in this rotating frame.  These dynamical
variables are given by the following expressions:
\be
K(t)=\frac{1}{2}\om^2 \sum_{j=1}^{j=n}(\frac{dx_j}{dt})^2+
(\frac{dy_j}{dt})^2+ (\frac{dz_j}{dt})^2,
\en
\be
\label{potential1}
P(t)=-\sum_{i,j}\frac{1}{r_{ij}(t) ^3},
\en
where $r_{ij}(t)=\sqrt{(x_i(t)-x_j(t))^2+(y_i(t)-y_j(t))^2 +(z_i(t)-z_j(t))^2}$
is the distance between particles  $i$ and $j$,
\be
\label{iner1}
I_x(t)=\sum_{i=1}^{i=n} y_(t)_i^2+z_i(t)^2,
\en
is the moment of inertia, and
\be
\label{angm}
L_x(t)=\om \sum_{i=1}^{i=n} y_i(t) \frac{dz_i(t)}{dt}-z_i(t) \frac{dy_i(t)}{dt}
\en
is the angular momentum.
The extrema of the action integral, Eq. \ref{action1}, is obtained by
determining the Fourier coefficients for which the partial 
derivatives of the action 
with respect to these coefficients vanish \cite{michael1},

\be
\label{ao}
\frac{\partial A}{\partial a(k_o)} = u(k_0)a(k_o)-f(k_o)=0,
\en

\be
\label{be}
\frac{\partial A}{\partial b(k_e)} = u(k_e)b(k_e)-g(k_e)+ v(k_e)c(k_e)=0,
\en
                                                                                           
\be
\label{ae}
\frac{\partial A}{\partial c(k_e)} = u(k_e)c(k_e)-h(k_e)+ v(k_e)b(k_e)=0,
\en
                                                                                           
where $u(k)=k^2\om^2+\Om^2$ and  $v(k)=2k\om \Om $, and 
\be
\label{fso}
f_(k)=\frac{1}{\pi}\int_0^{2\pi}dt
\frac{\partial P}{\partial x}sin(kt),
\en

\be
\label{gse}
g(k)=\frac{1}{\pi}\int_0^{2\pi}dt
\frac{\partial P}{\partial y}sin(kt),
\en

\be
\label{gco}
h(k)=\frac{1}{\pi}\int_0^{2\pi}dt
\frac{\partial P}{\partial z}cos(kt).
\en

These equations are the Fourier transform
of the differential  equations  of motion, but 
for our steepest descent gradient method of solution
it is essential  to recognize   
that these equations determine  also an
extrema of the action integral, Eq. \ref{action1}, even when  
only a {\it finite} number of Fourier terms in the expansion
the the orbital coordinates are taken into account. 
For convergent Fourier expansions, the number of such terms 
depends, of course, on the desired numerical accuracy for the solutions.
Starting with the Fourier coefficients to some  approximate 
form of the periodic  orbit, and evaluating the  partial derivatives
of the action, Eqs. \ref{ao} -\ref{ae},
an improved  orbit is then  obtained by changing
the value of each  coefficient in  proportion to 
the  corresponding partial derivative of the action.
For example,  with some suitable initial values for 
the coefficient $a(k_o)$, $b(k_e)$ and $c(k_e)$, we obtain a new value $a'(k_o)$
by the transformation 
\be
\label{map1}
a'(k_o)=a(k_o)+\frac{\del s}{u(k_o)} \frac{\partial A}{\partial a(k_o)},
\en
with similar expressions for the other Fourier coefficients $b(k_e)$ and $c(k_e)$,
where $\del s$ is a parameter. Then, to first order in the difference
$\del a(k_o)=a'(k_o)-a(k_o)$, the change $\del A$ in the action integral is given by
\be 
\del A= \del a(k_0)    \frac{\partial A}{\partial a(k_o)}=
\frac{\del s}{u(k_o)}(\frac{\partial A}{\partial a(k_o)})^2.
\en
The sign of the parameter $\del s$ must be chosen to
be positive (negative) depending on whether  the action 
is a maximum (minimum) with respect to the variation $\del a(k_o)$. 
This procedure is then iterated until the partial derivatives
of the action are reduced to any  desired accuracy.
Of course, the sign of $\del s$ is not known a-priori, but
if the wrong value is chosen, one finds that the iteration
diverges. 
Our method differs in an important  way from the standard 
method of steepest descent by the introduction of the factor $1/u(k)$ 
in Eq. \ref {map1}, which is important to obtain  convergence
with moderate values of the paramater  $\del s$  
(in the range $0.1$ to $0.5$). In this range 
most of our calculations have converged after a 
a few hundred interations to values for 
$\del a=a'(k_o)-a(k_o)$ of order $10^{-10}$
to $10^{-13}$, 
and similar  for the other Fourier coefficients. 
However,  near the limit  $\Omega=1$, which  corresponds to
the Lagrange circular orbit, we found that 
about $10^4$ to $10^5$ iterations were necessary 
to obtain  convergence of the dominant Fourier  coefficient
$a(1)$ to only  few decimal places. The reason for this 
very slow convergence turns out be 
that in the  Lagrange limit, the transformation  given by  Eq. \ref{map1} 
has a marginal eigenvalue, which
will be discussed in the next section. For a frame 
rotating around the $x$ axis of symmetry, we found that
the various extrema of the action integral always occur at a
{\it minimum} of this integral for variations 
with respect to the Fourier coefficients shown in Figs. 4-6, 
but for rotations around the $y$ axis, the extrema  with respect to the
coefficients  $c(1)$ and $c(3)$ occurred at a {\it maximum} 
of the action integral. For rotations around the
$y$ axis of symmetry of the figure eight (see Figs. 7 and 8), 
we found  that for  $\Omega$ somewhat greater than $0.8$, 
the iterations of our Fourier gradient 
transformation  diverged,  suggesting the existence
of a nearby singularity in the dependence of the
corresponding family of solutions as a function of $\Omega$.   

\subsection *{II. Perturbation Solution Near the Lagrange Circular Orbit}

To illustrate  the application of our steepest descent Fourier gradient 
transformation, Eq.\ref{map1},
we calculate analytically the integrals in Eqs. \ref{fso}-\ref{gco}
near the Lagrange solution for $\Omega=1$.
Setting $\eps =a(1)/r $, $b_2=b(2)/r $, 
and $b_4=b(4)/r$  where $r=1/3^{1/6}$ is the length scale
for $m=1$, we  evaluate these integrals 
by perturbation theory to first order in $\beta=1-\Omega$, and obtain 
the following transformations
for the coefficients $\eps,b_2$ and $b_4$:
\be
\eps'=\frac{\eps}{ b_2^3}(1-\frac{9}{8}\eps^2+\frac{3}{2}b_4),
\en
\be
b_2'=(4(1-\beta)b_2+\frac{1}{ b_2^2}(1-\frac{3}{4}\eps^2))/(5-2\beta),
\en
and
\be
b_4'=(8(1-\beta)b_4+\frac{1}{b_2^2}(\frac{3}{8}\eps^2-b_4))/(17-2\beta).
\en

The fixed points  for these transformations are given by 
\be
\eps_ f=\sqrt{\frac{19}{3}\beta},
\en
\be
b_{2f}=1-\frac{9}{4}\beta,
\en
and
\be
b_{4f}=\frac{1}{4}\beta.
\en
These fixed points correspond to the lowest order  
perturbation solution near the Lagrange circular orbit,
previously obtained by C. Marchal \cite{march1}, \cite{march2}, 
except for the sign for $b_{4f}$, which is given incorrectly 
in these references.

It is illuminating to calculate the first order deviations 
$\del \eps =\eps -\eps_f$, $\del b_2=b_2'-b_{2f}$
and $\del b_4=b_4'-b_{4f}$, which satisfy 
the simpler equations
\be
\del \eps'=\del \eps(1-4\beta -6\sqrt{3\beta/19}  \del \eps ),
\en
\be
\del b_2'=\frac{2}{5}(1-\frac{42}{5}\beta)\del b_2,
\en
\be
\del b_4'=\frac{15}{34}(1-\frac{274}{255}\beta)\del b_4.
\en
Evidently, in the Lagrange limit where $\beta$ vanishes,  
iterations of $\del b_2$ and $\del b_4$ converge rapidly
to the respective fixed points,
while the transformation  for $\eps$ has an
eigenvalue which approaches unity, and therefore converges increasingly
slowly as $beta$ decreases.  
This property explains the slow convergence which we found for iterations
of the exact transformations   near the Lagrange limit. For example, 
if one wishes to determine the fixed point to an
accuracy $\del =10^{-m}$, then the number of iterations 
required for  this map is $n \approx 2.3*m/4\beta$. 
For example, to obtain $\eps_f$ for $\beta=.0001$  to an accuracy 
$\del=10^{-4}$,  which gives  
$\eps_f=.025166$,  one needs approximately $n=23,000$ iterations 
of these linearized equations.

\subsection *{ III. Numerical results}

In this section we present our numerical results in graphical form  for 
finite-angular momentum figure eight shaped orbits 
in a rotating frame with angular frequency $\Omega$.
These orbits are obtained by the steepest descent gradient 
Fourier method described
in section I. For rotations along the horizontal symmetry axis ( $x$-axis) 
of the planar figure eight orbit, the values of  $\Omega$ 
are in the range $0$ to $1$.
As an illustration of the characteristic shape of such orbits
in the rotating frame, we show in Fig.1 a three dimensional  finite angular
momentum  periodic orbit for three equal masses in a frame rotating 
with frequency $\Omega=0.5$. 
The shape of this orbit corresponds to a  three
dimensional figure eight which is bend symmetrically downward  at a 
self-crossing point
of the orbit that lies on the $z$ axis  at $x=y=0$,
as can be seen from Eqs.\ref{xt}-\ref{zte}.
The position of each of the three masses  is shown at five equal
intervals during  a third of the period $\tau=2 \pi$, in corresponding
segments of the orbit indicated by blue, red and green colors. The segments
of the orbit below the $x-y$ plane are indicated by dashed lines.
For positive values of $x$ this orbit is transversed clockwise as
viewed on a projection on the $x-y$ plane,  while for
negative values of $x$ it is transversed counterclockwise, where
the entire  frame is rotating  counterclockwise around the x-axis
with respect an  inertial frame.

We have obtained similar orbits for an odd number of equal masses $n\geq 3$. 
A typical case for a finite angular momentum orbit 
with $21$ masses is shown in red  in Fig. 2, 
 while the associated  planar figure eight orbit for 
$\Omega=0$ is shown in blue, 
and the Lagrange circular orbit for $\Omega=1.0$ is shown in green. 
In Fig. 3 we show the projection on the $x-y$ plane 
of six of these orbits, with  $n=3$ and $\Omega = 
0, 0.2, 0.4, 0.6, 0.8$ and $1.0$,
corresponding to similar figures in \cite{mont2}.
In Figs. 4-6, we plot the first few Fourier coefficients of the
expansion of the coordinates for the finite
angular momentum figure eight solution, Eqs. \ref{xt}-\ref{zto}, 
as a function of 
the angular frequency $\Omega$, which give  numerical evidence 
that these coefficients are continuous functions 
of $\Omega$ between $\Omega=0$, for the
coefficients for the planar figure eight orbit, and $\Omega=1$,
for the coefficient of the classical Lagrange orbit.
We show with dashed lines Marchal's perturbation approximation 
for these coefficients \cite{march1},\cite{march2} 
near the Lagrange limit $\Omega=1$, evaluated  to second order 
in powers of the square root of $\beta=1-\Omega$, and also our
corresponding expansion to second order in
powers of $\Omega$ near  $\Omega=0$  for the
bifurcation near the planar figure eight orbit \cite{michael1}.

We have previously evaluated also some orbits for 
finite angular momentum  of the figure
eight orbit rotating about its other two symmetry 
axes \cite {michael1}, and  
recently existence proofs and  further  numerical 
results \cite{mont2} have extended this work.
In Fig.7 we show  an orbit for rotation about
the $y$ axis with $\Omega=0.5$. In this case the
orbit does not cross at $x=y=0$; instead the
two values of $z$ are equal but opposite in sign
as can be readily seen from Eqs. \ref{xt},\ref{yt} and \ref{zto}.
In Fig. 8 we show  the projection of
these orbits  on the $x-y$ plane 
for $\Omega=0., 0.2, 0.4, 0.6$ and $0.8$,
corresponding to similar results in \cite{mont2}.
This result indicates  that
in this case the projected  crossing angle 
at $x=y$ is a constant.
But for values of $\Omega$ greater than about 0.8 we found that 
our iteration procedure  ceases to
converge, suggesting that
there is a singularity in the dependence of these
coefficient on $\Omega$.
Finally, we give in Fig.9  an illustration of the shape of the
orbit when the axis of rotation is normal to the
plane of the figure eight orbit, similar to one
shown in \cite{michael1}. In this case the finite
angular momentum orbit remains planar but the symmetry 
about the $y$ axis is broken.
 
We  remark that our steepest descent gradient method also converges   
when the axis of the rotating frame is not taken  along any of the
symmetry axes of the planar figure eight orbit, 
but this procedure  does not generate new solutions.
In this case one has to take into account
additional terms in the Fourier expansion of the coefficients 
of the orbit, and the previous  separation  into even
and odd integers in the sums, Eqs. \ref{xt}-\ref{zto}  is no longer valid.
We find that  for not too small values of $\Omega$, 
the iterations converge to a solution
which correspond to a  member of the family associated with  
rotations along the horizontal symmetry axis of the figure eight orbit,
which is now located along the chosen rotation axis.
But as  $\Omega$ approaches $0$ the convergence 
of the iterations becomes extraordinarily
slow; e.g. for $\Omega=.01$ one requires about $10^6$ iterations
to achieve convergence. 

\subsection*{IV Nonlinear Stability}

In this section we present numerical results for the
stability of the three families of periodic figure eight orbits 
with finite angular momentum along the
three symmetry axis of the eight. We also give  some examples 
of stable orbits located in the vicinity of these orbits
near the critical point where the periodic orbits of a 
family become unstable. To determine
the stability  of a periodic  orbit we
calculate it numerically with 
the fourth order Runge-Kutta algorithm, taking as initial  
coordinates and momenta  values very close  to those  obtained 
from our Fourier gradient method (see section I). 
We then follow this orbit for up to approximately $ 10^3$
periods monitoring the constancy of the total  energy and the 
total angular moment to an accuracy  $10^{-5}$.  
For the family of orbits rotating around the horizontal
axis of symmetry of the eight, a linear stability analysis 
in reference \cite{mont2},  reported that these orbits
are linearly unstable for all finite angular frequencies $\Omega>0$. 
We find, however, that for $\Omega$ approximately less than
0.1 the orbits bifurcate
into {\it stable} albeit  quite complex quasi-periodic orbits.
An example for such a stable nearby orbit for $\Omega=.05$
is shown in Fig. 10,  which is plotted for $60$ periods.
Each of the the three equal masses  moves  on
separate but nearby orbits indicated by the colors blue, red
and green. The separation between these orbits
decreases monotonically as $\Omega$ decrease.
In Fig. 11 the corresponding orbit is shown after $120$
periods. As the number of periods increases, further wings appear
and the orbit starts to fill a spherical volume with 
the radius of half the horizontal axis of the figure eight
orbit.  Evidently, in this case the linear 
stability analysis  given in reference \cite{mont2} is inadequate to determine
 the domain for the existence of stable orbits in the
vicinity of this family. But for the family 
of periodic figure eight orbit  rotating  around the $y$ and the $z$ axis, 
our numerical stability results are in accordance with the 
the linear analysis \cite{mont2}. As an example, in Fig. 12 
we show a quasi-periodic period orbit for $\Omega_y=.09$ for $166$
periods. We have checked  for $333$ periods that 
for   $\Omega_y=.0920$ such 
an orbit remain stable, while for larger values of $\Omega_y$  
it becomes unstable.
Likewise, for rotations around the $z$ axis of symmetry
of the eight we found  nearby quasi-periodic orbits for 
$\Omega_z <.585$ which we verified to remain stable
for $333$ periods.  For $\Omega_y=.55$ an example 
is shown in Fig. 13 for $20$ periods. 
These orbits correspond to  a slow precesion with small deformations
of the stable orbit shown in Fig. 10.

\subsection*{Acknowledgements}
I would like to thank C. Marchal for his comments.

\pagebreak

\begin{figure}
\begin{center}
\epsfxsize=\columnwidth
\epsfig{file=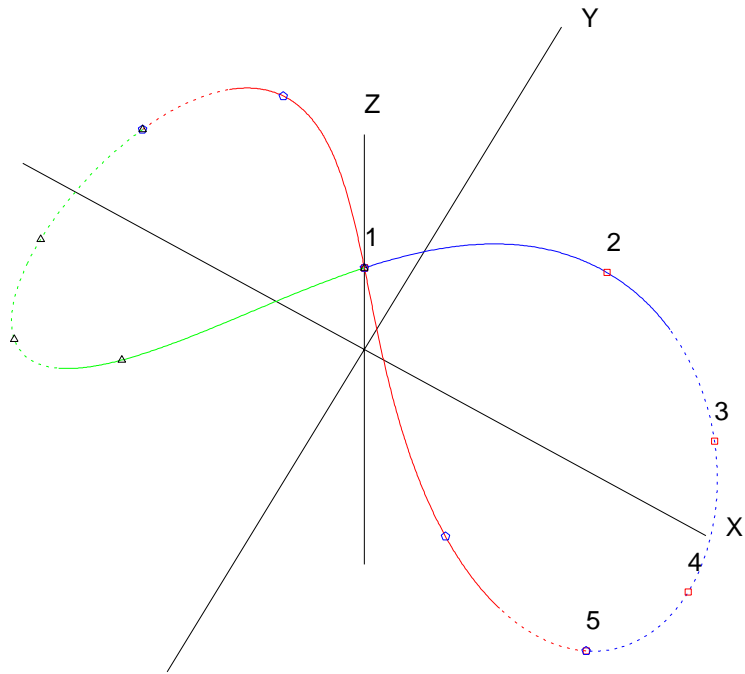, width=15cm}
\end{center}
\caption{ Finite angular momentum figure eight orbit, which is periodic 
in a frame rotating around the $x-axis$
with frequency $\Omega=0.5$. The position of the three masses are
indicated by squares at 4 equal time intervals during 1/3 of a period. 
}
\label{Fig. 1}
\end{figure}

\begin{figure}
\begin{center}
\epsfxsize=\columnwidth
\epsfig{file=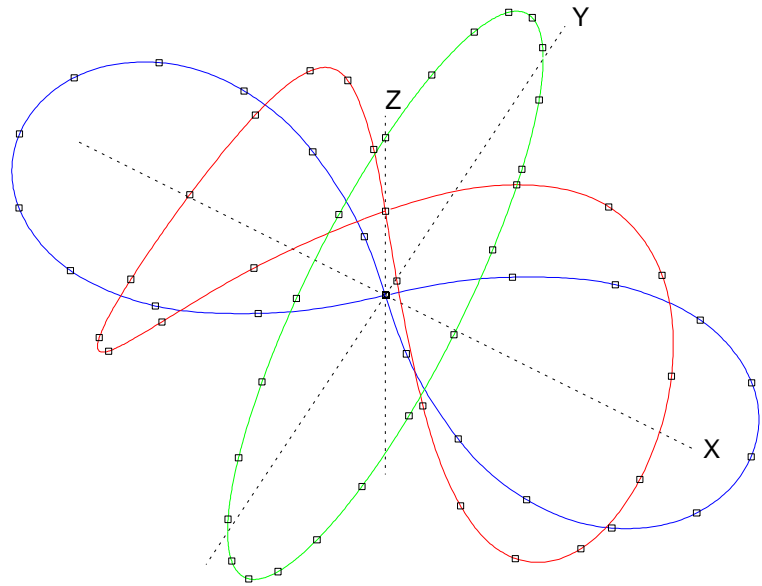, width=15cm}
\end{center}
\caption{ Periodic orbits with 21 particles of equal
mass shown at equal time intervals indicated
by squares. The blue curve is the planar figure eight orbit
in the $x-y$ plane, 
the green curve is the circular Lagrange orbit in the
$y-z$ plane , and the
red orbit is a member of a family of 
finite angular momentum figure 
eight orbits which is periodic in a frame rotating with
frequency $\Omega=0.5$ around the x-axis.
}
\label{Fig. 2}
\end{figure}

\begin{figure}
\begin{center}
\epsfxsize=\columnwidth
\epsfig{file=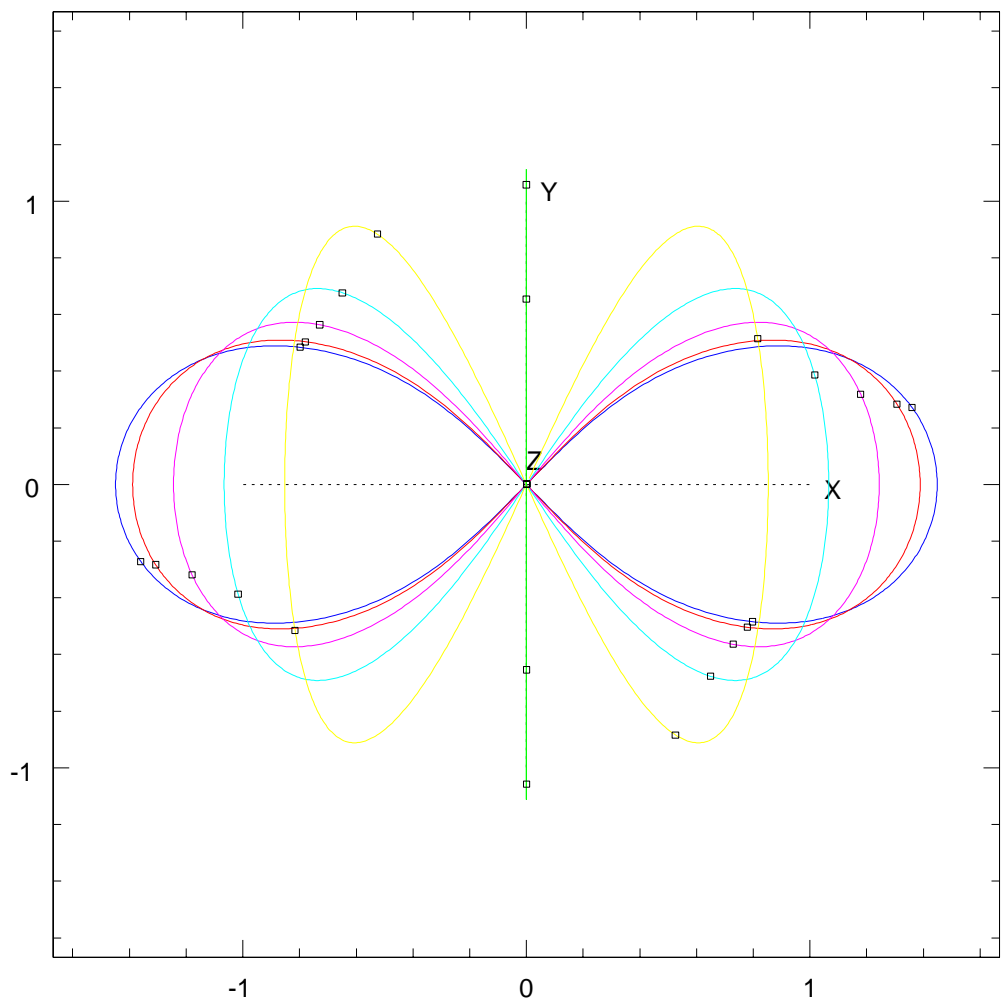, width=15cm}
\end{center}
\caption{ Projection on the $x-y$ plane of finite angular momentum figure
eight orbits for three equal masses rotating shown in a frame around
the x-axis of symmetry for $\Omega= 0,0.2,0.4,0.6,0.8$ and $1.0$
(from blue to green).
}
\label{Fig. 3}
\end{figure}

\begin{figure}
\begin{center}
\epsfxsize=\columnwidth
\epsfig{file=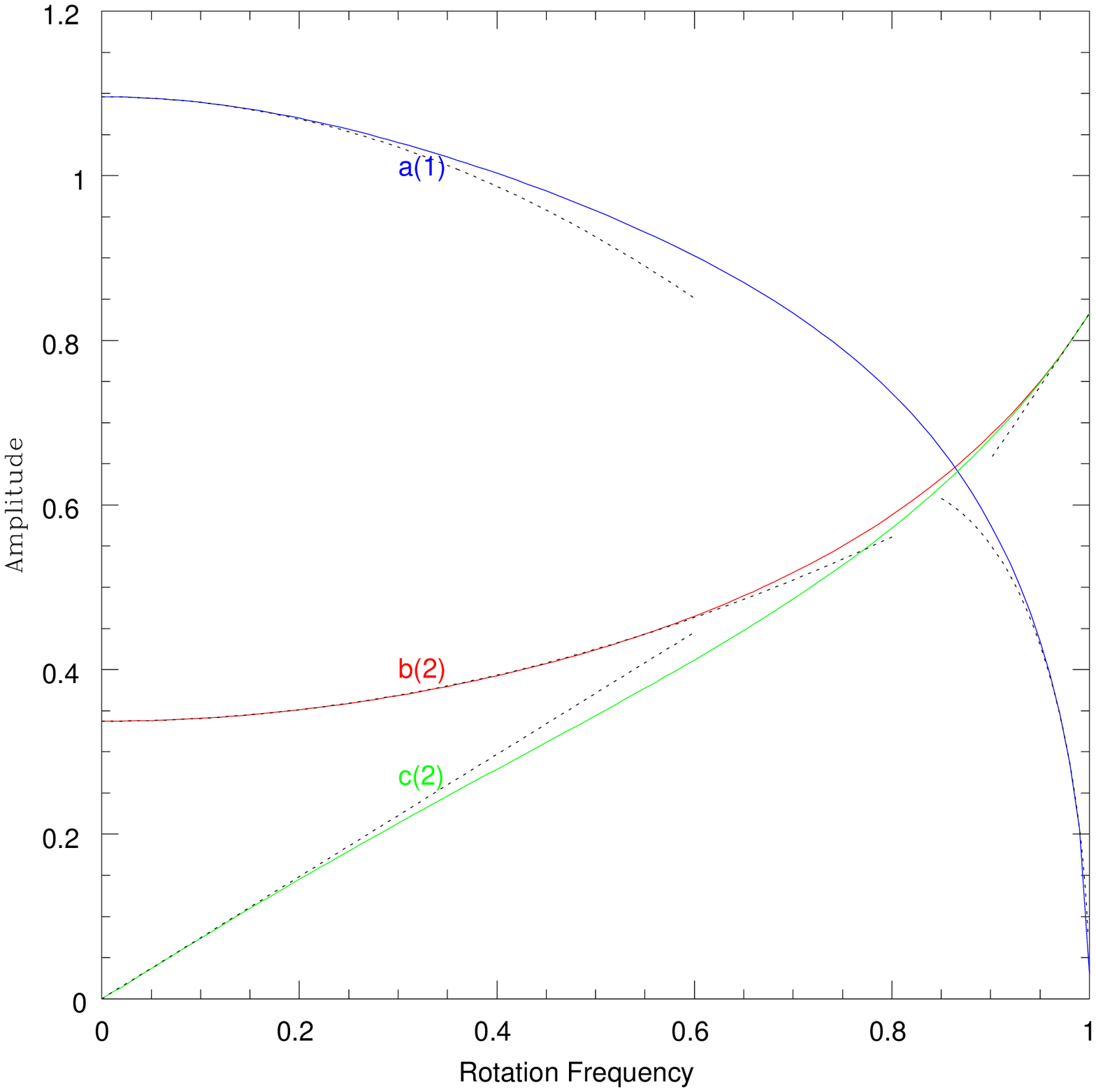, width=15cm}
\end{center}
\caption{ First three Fourier coefficients for the finite angular momentum
figure eight solutions  as functions of the rotation frequency
$\Omega$. The dashed curves are perturbation solutions.
}
\label{Fig. 4}
\end{figure}

\begin{figure}
\begin{center}
\epsfxsize=\columnwidth
\epsfig{file=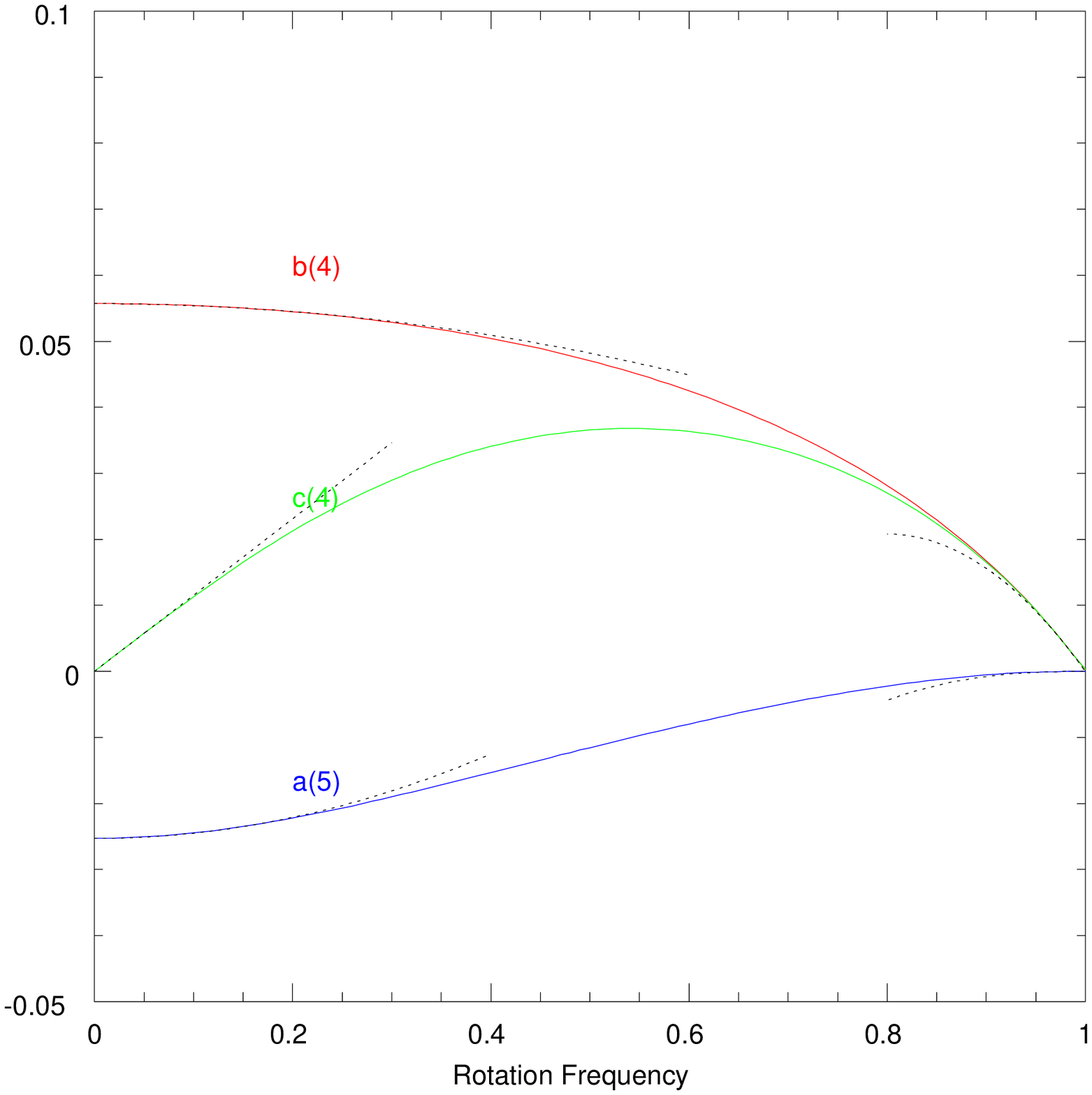, width=15cm}
\end{center}
\caption{
Fourier series  coefficients for the finite angular momentum
figure eight solutions  as functions of the rotation frequency
$\Omega$. The dashed curves are perturbation solutions (see text)
}
\label{Fig. 5}
\end{figure}

\begin{figure}
\begin{center}
\epsfxsize=\columnwidth
\epsfig{file=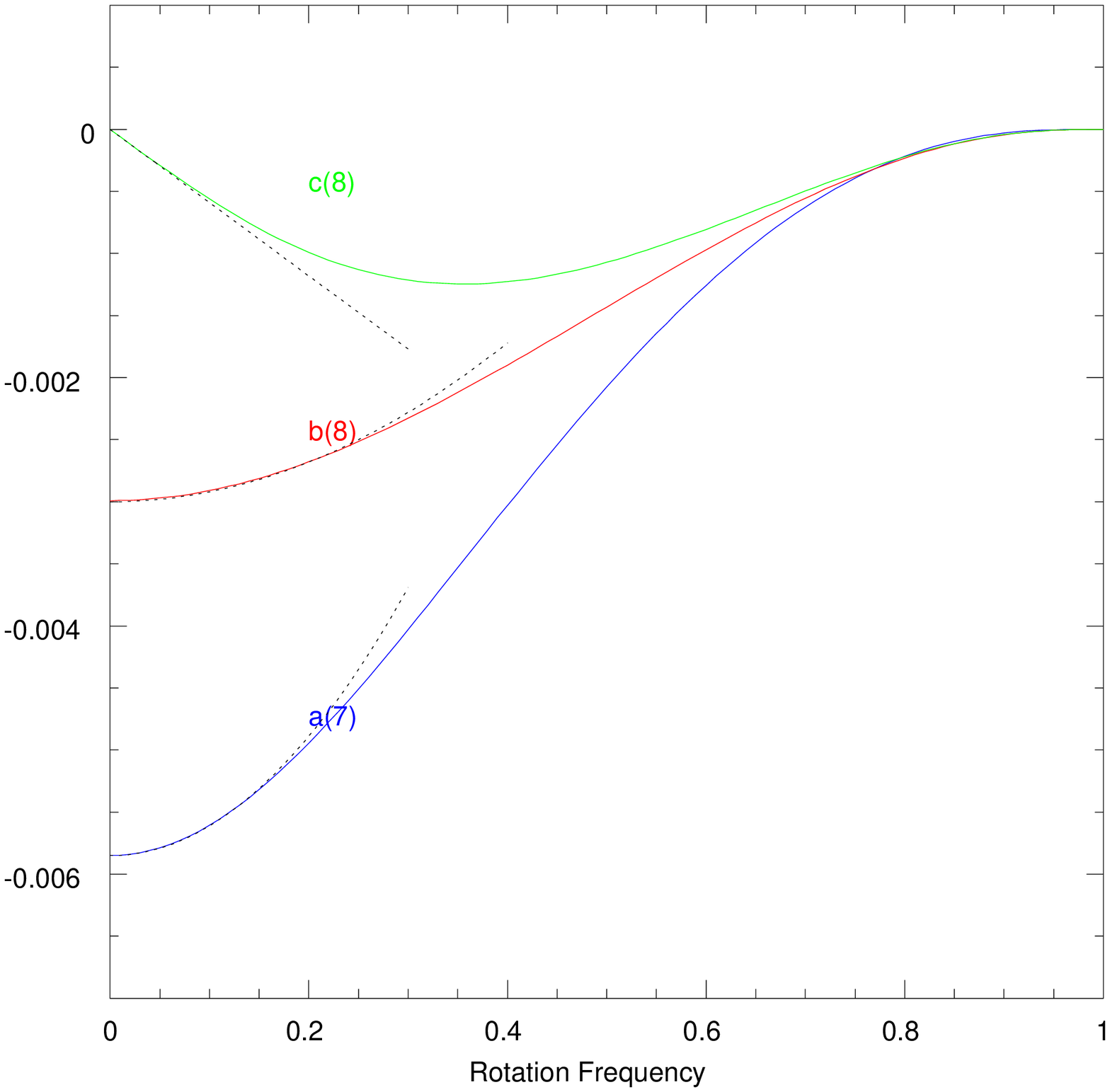, width=15cm}
\end{center}
\caption{
Fourier series  coefficients for the finite angular momentum
figure eight solutions  as functions of the rotation frequency
$\Omega$. The dashed curves are perturbation solutions.
}
\label{Fig. 6}
\end{figure}

\begin{figure}
\begin{center}
\epsfxsize=\columnwidth
\epsfig{file=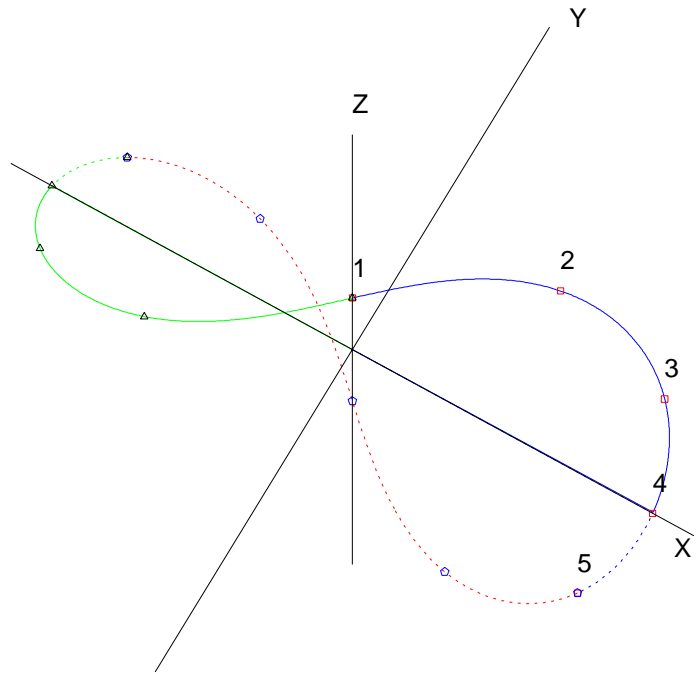, width=15cm}
\end{center}
\caption{ Finite angular momentum figure eight orbit rotating 
around the $y$ axis of symmetry for $\Omega=0.5$. 
The position of the three masses are indicated 
by squares at 4 equal time intervals during 1/3 of a period. 
}
\label{Fig. 7}
\end{figure}

\begin{figure}
\begin{center}
\epsfxsize=\columnwidth
\epsfig{file=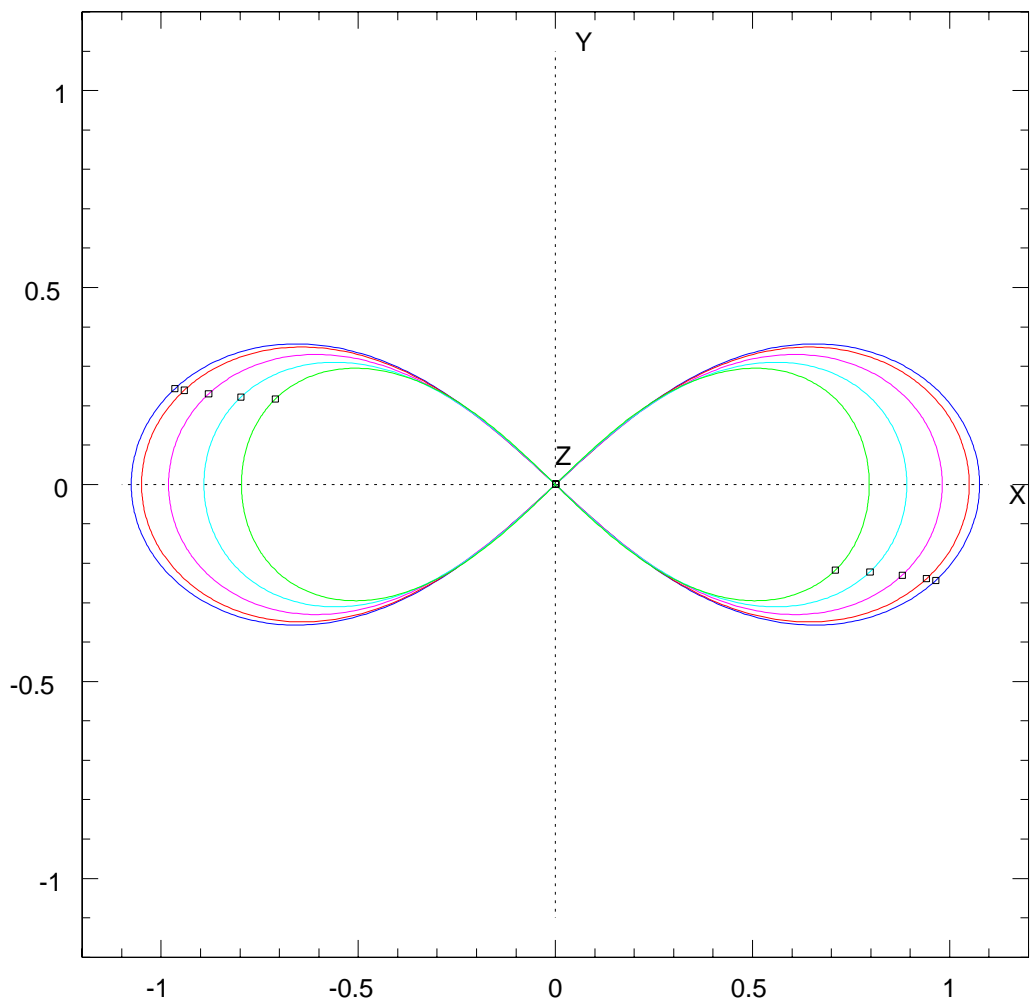, width=15cm}
\end{center}
\caption{ Projection on the $x-y$ plane of finite
angular momentum figure eight solutions shown in a frame rotating 
around the $y$ axis for $\Omega=0,0.2,0.4,0.6,$ and $0.8$
(from blue to green).
}
\label{Fig. 8}
\end{figure}

\begin{figure}
\begin{center}
\epsfxsize=\columnwidth
\epsfig{file=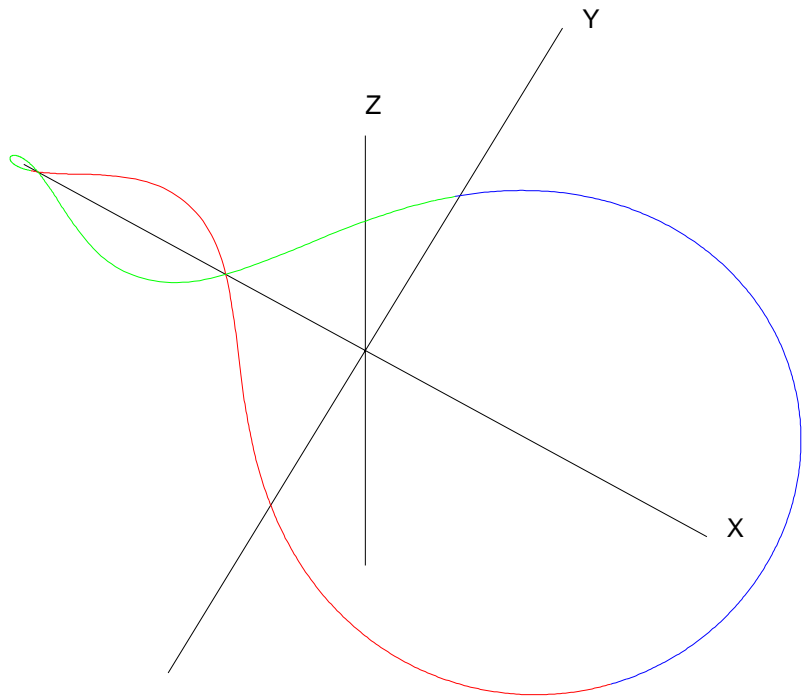, width=15cm}
\end{center}
\caption{Figure eight orbit  rotating around the $z$ axis with
$\Omega=.5$. The  blue,red and green segments are the trajectories
for  three equal masses during  one third of a period. 
}
\label{Fig. 9}
\end{figure}

\begin{figure}
\begin{center}
\epsfxsize=\columnwidth
\epsfig{file=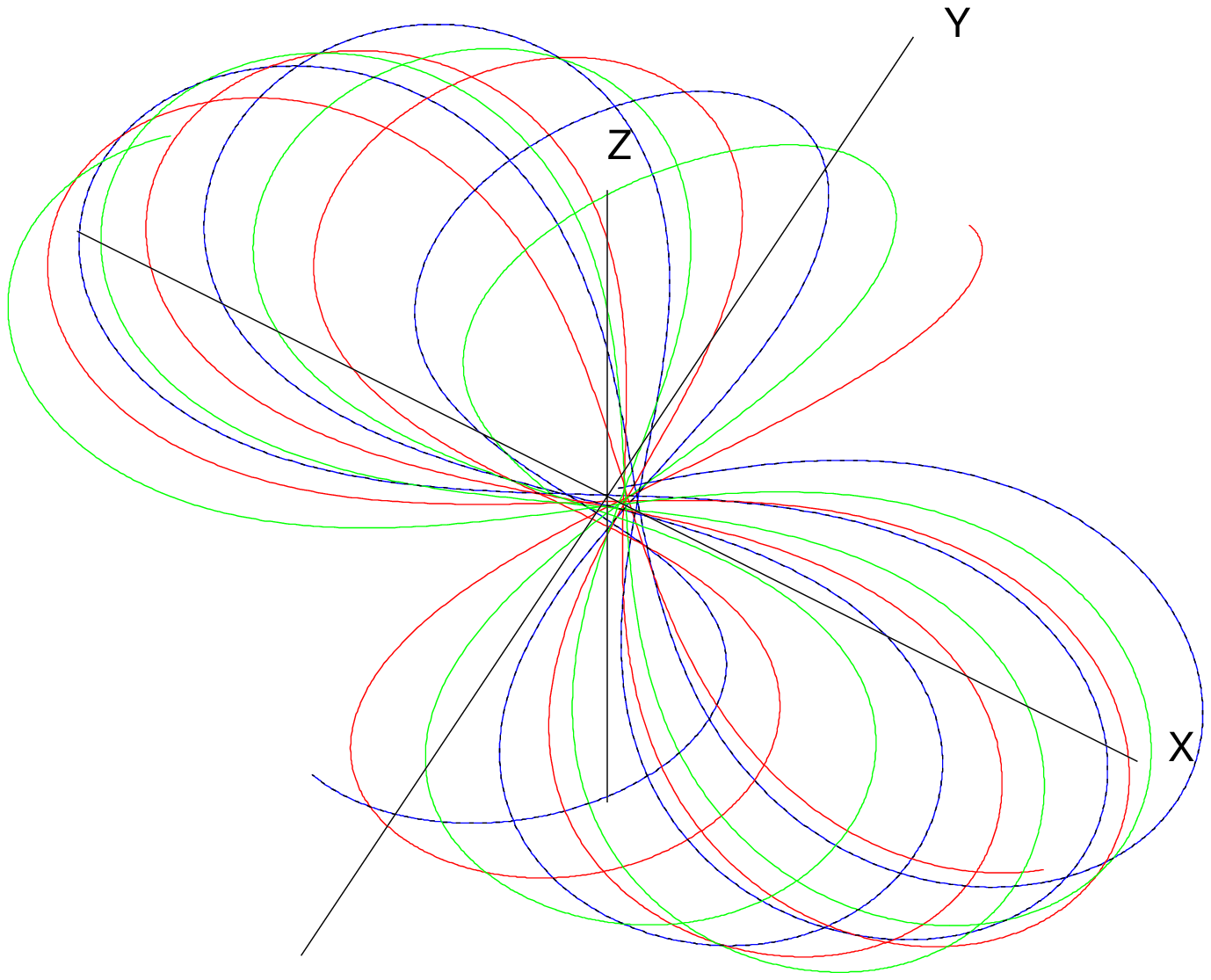, width=15cm}
\end{center}
\caption{Stable figure  eight orbit  in a rotating frame around the 
$x$ axis with $\Omega=.05$. The  blue,red and green  orbits
are the separate trajectories for  three equal masses. 
}
\label{Fig. 10}
\end{figure}

\begin{figure}
\begin{center}
\epsfxsize=\columnwidth
\epsfig{file=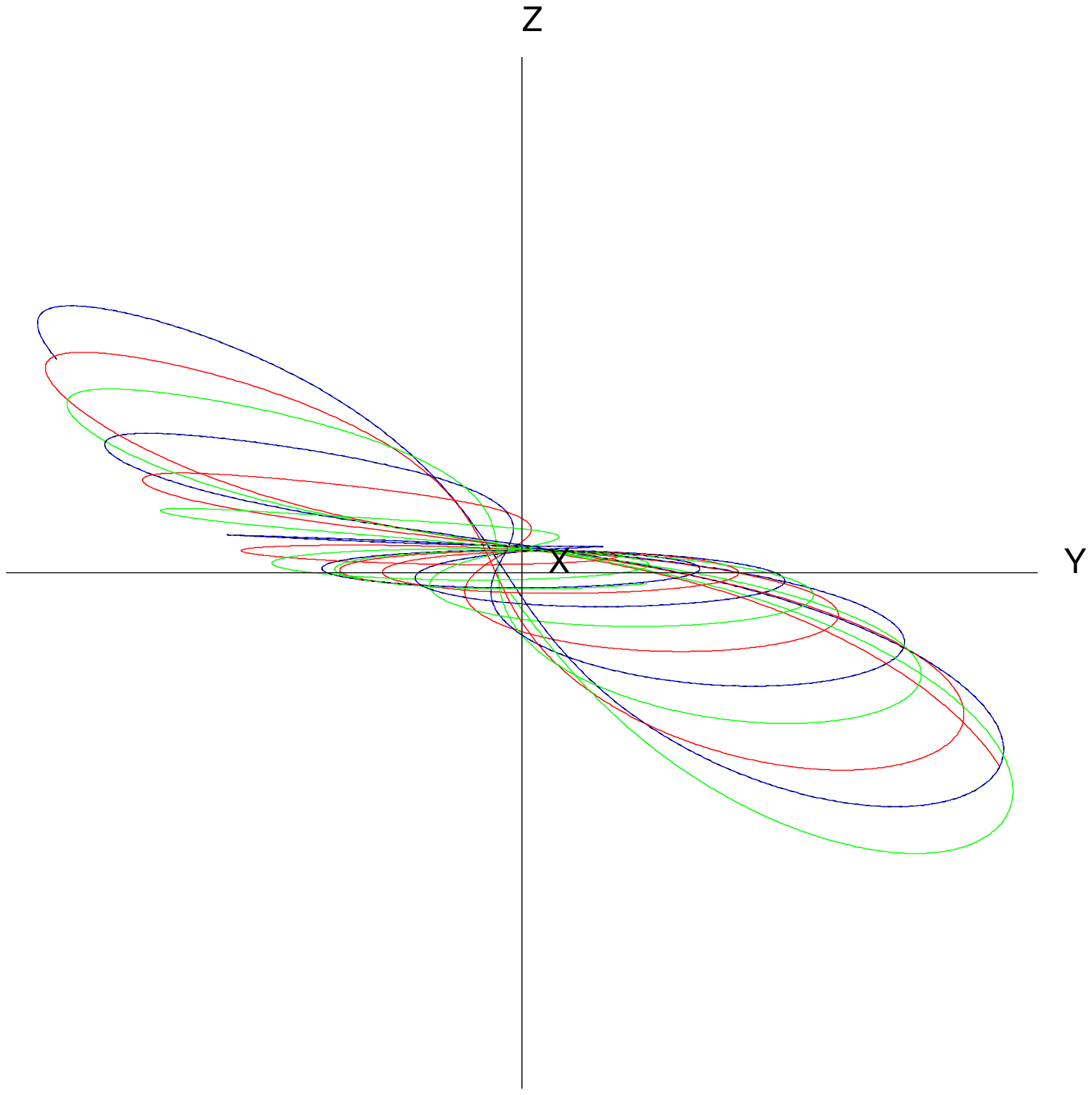, width=15cm}
\end{center}
\caption{Stable figure  eight orbit  in a rotating frame around the 
$x$ axis with $\Omega=.05$. The  blue,red and green  orbits
are the separate trajectories for  three equal masses. 
}
\label{Fig. 11}
\end{figure}

\begin{figure}
\begin{center}
\epsfxsize=\columnwidth
\epsfig{file=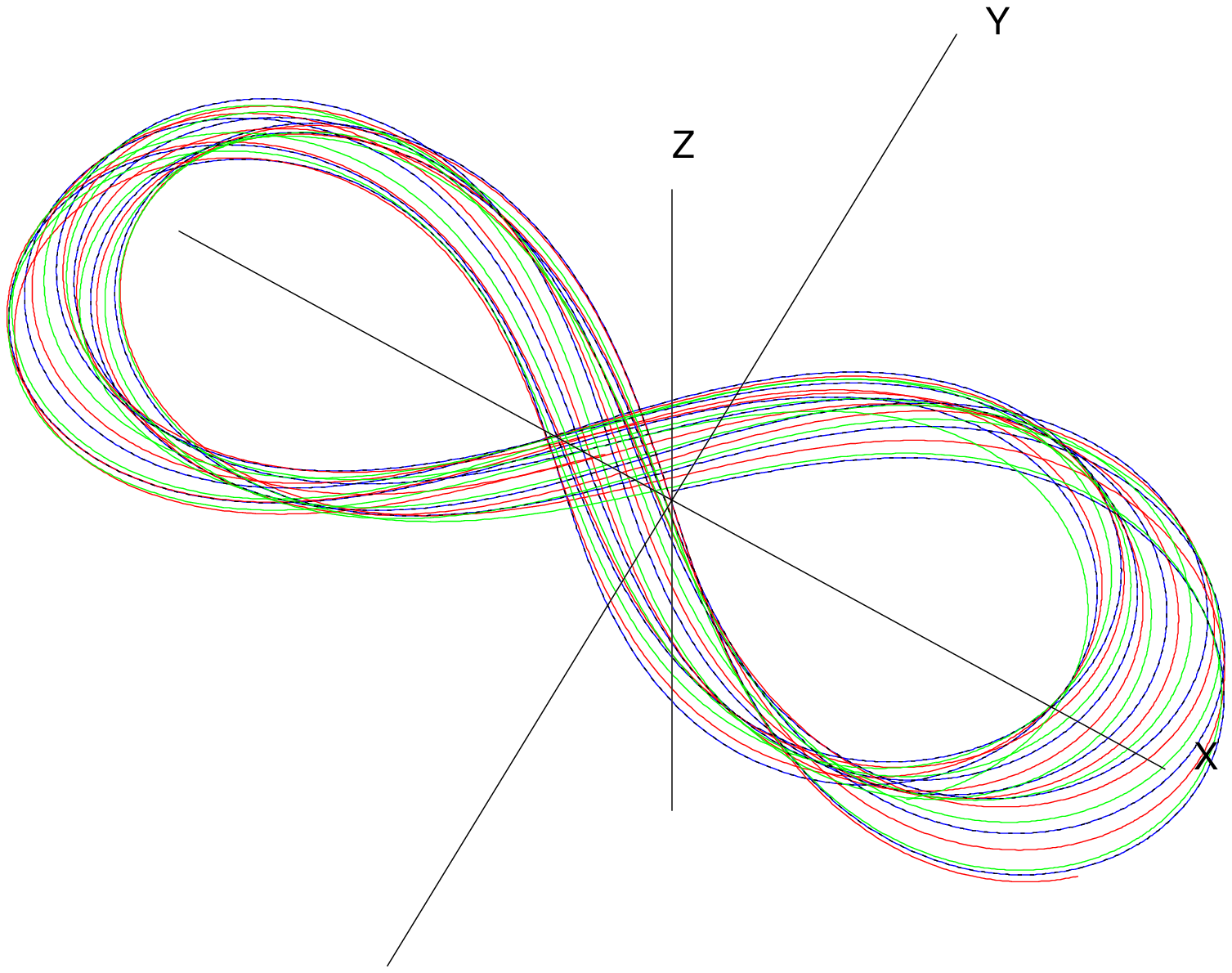, width=15cm}
\end{center}
\caption{Stable figure eight orbit  in a rotating frame around the 
$y$ axis with $\Omega=.05$. The  blue,red and green  orbits
are the separate trajectories for  three equal masses.
}
\label{Fig. 12}
\end{figure}

\begin{figure}
\begin{center}
\epsfxsize=\columnwidth
\epsfig{file=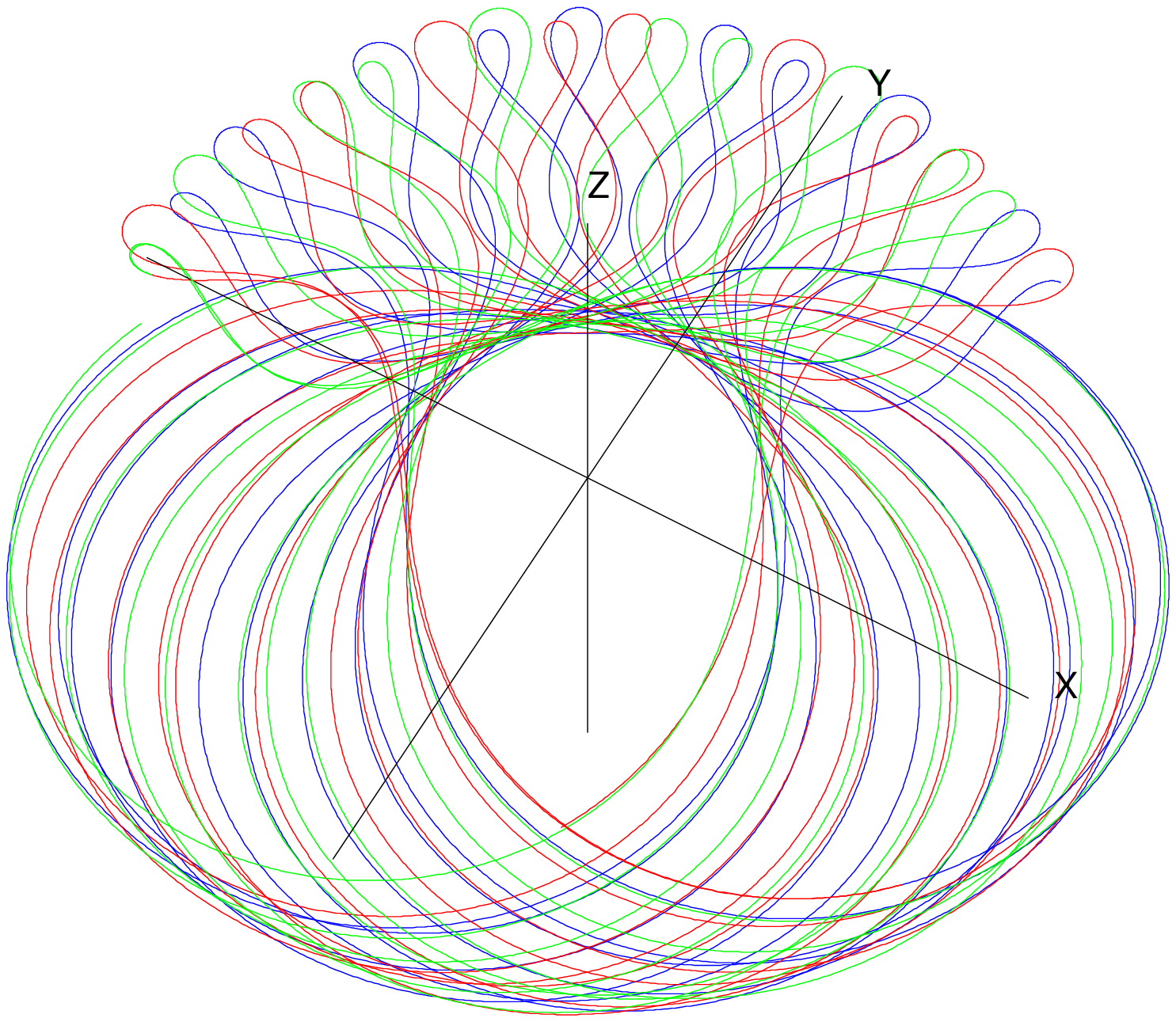, width=15cm}
\end{center}
\caption{Stable figure eight orbit  in a rotating frame around the 
$z$ axis with $\Omega=.55$. The  blue,red and green  orbits
are the separate trajectories for  three equal masses.
}
\label{Fig. 13}
\end{figure}
\end{document}